\input amstex
\input amsppt.sty
\magnification=\magstep1
%\textwidth=13.5cm
%\textheight=24cm
%\hoffset=-1cm
%\hsize=30truecc
%\baselineskip=16truept
%\vsize=22.2truecm
\hsize=33truecc
\vsize=22.2truecm
\baselineskip=16truept
\NoBlackBoxes
\nologo
\pageno=1
\topmatter
\TagsOnRight

\def\N{\Bbb N}
\def\Z{\Bbb Z}

\def\b{\bigg}

\def\({\b(}
\def\[{\b[}
\def\){\b)}
\def\]{\b]}

\def\t{\text}

\def\em{\emptyset}
\def\se {\subseteq}

\def\sm{\setminus}

\def\eq{\equiv}

\def\ls{\leqslant}
\def\gs{\geqslant}

\def\ve{\varepsilon}
\def\da{\delta}

\def\Proof{\noindent{\it Proof}}
\def\Remark{\noindent{\it Remark}}

\def\Ack{\noindent {\bf Acknowledgments}}
\hbox{J. Number Theory 171(2017), 275--283.}
\medskip
\title On $x(ax+1)+y(by+1)+z(cz+1)$ and $x(ax+b)+y(ay+c)+z(az+d)$\endtitle
\author Zhi-Wei Sun \endauthor
\affil Department of Mathematics, Nanjing University
     \\Nanjing 210093, People's Republic of China
    \\  zwsun\@nju.edu.cn
    \\ {\tt http://math.nju.edu.cn/$\sim$zwsun}
 \endaffil
\abstract In this paper we first investigate for what positive integers $a,b,c$ every nonnegative integer $n$ can be written as $x(ax+1)+y(by+1)+z(cz+1)$ with $x,y,z$ integers.
We show that $(a,b,c)$ can be either of the following seven triples
$$(1,2,3),\ (1,2,4),\ (1,2,5),\ (2,2,4),\ (2,2,5),\ (2,3,3),\ (2,3,4),$$
and conjecture that any triple $(a,b,c)$ among
$$(2,2,6),\ (2,3,5),\ (2,3,7),\ (2,3,8),\ (2,3,9),\ (2,3,10)$$
also has the desired property. For integers $0\ls b\ls c\ls d\ls a$ with $a>2$, we prove that any nonnegative integer
can be written as $x(ax+b)+y(ay+c)+z(az+d)$ with $x,y,z$ integers, if and only if the quadruple $(a,b,c,d)$ is among
$$(3,0,1,2),\ (3,1,1,2),\ (3,1,2,2),\ (3,1,2,3),\ (4,1,2,3).$$
\endabstract
\thanks 2010 {\it Mathematics Subject Classification}.
Primary 11E25; Secondary 11D85, 11E20.
\newline \indent Supported by the National Natural Science
Foundation (Grant No. 11571162) of China.
\endthanks
\endtopmatter
\document

\heading{1. Introduction}\endheading

Let $\N=\{0,1,2,\ldots\}$. Throughout this paper, for $f(x,y,z)\in\Z[x,y,z]$ we set
$$E(f(x,y,z))=\{n\in\N:\ n\not=f(x,y,z)\ \t{for any}\ x,y,z\in\Z\}.$$
If $E(f(x,y,z))=\em$, then we call $f(x,y,z)$ {\it universal over $\Z$}.
The classical Gauss-Legendre theorem (cf. [N96, pp.\,3-35]) states that
$$E(x^2+y^2+z^2)=\{4^k(8l+7):\ k,l\in\N\}.$$

Recall that those $T_x=x(x+1)/2$ with $x\in\Z$ are called {\it triangular numbers}.
As $T_{-x-1}=T_x$, $T_{2x}=x(2x+1)$ and $T_{2x-1}=x(2x-1)$, we see that
$$\{T_x:\ x\in\Z\}=\{T_x:\ x\in\N\}=\{x(2x+1):\ x\in\Z\}.\tag1.1$$

By the Gauss-Legendre theorem, any $n\in\N$ can be written as the sum of three triangular numbers
(equivalently, $8n+3$ is the sum of three odd squares). In view of (1.1), this says that
$$\{x(2x+1)+y(2y+1)+z(2z+1):\ x,y,z\in\Z\}=\N.\tag1.2$$
Motivated by this, we are interested in finding all those $a,b,c\in\Z^+=\{1,2,3,\ldots\}$
satisfying
$$\{x(ax+1)+y(by+1)+z(cz+1):\ x,y,z\in\Z\}=\N.\tag1.3$$

In the following theorem we determine all possible candidates $a,b,c\in\Z^+$ with (1.3) valid.

\proclaim{Theorem 1.1} Let $a,b,c\in\Z^+$ with $a\ls b\ls c$. If $x(ax+1)+y(by+1)+z(cz+1)$ is  universal over $\Z$,
then $(a,b,c)$ is among the following $17$ triples:
$$\aligned&(1,1,2),\,(1,2,2),\,(1,2,3),\,(1,2,4),\,(1,2,5),
\\&(2,2,2),\,(2,2,3),\,(2,2,4),\,(2,2,5),\,(2,2,6),
\\&(2,3,3),\,(2,3,4),\,(2,3,5),\,(2,3,7),\,(2,3,8),\,(2,3,9),\,(2,3,10).
\endaligned\tag1.4$$
\endproclaim
\Remark\ 1.1. As proved by Liouville (cf. [D99, p.\,23]),
$$\{2T_x+2T_y+T_z:\ x,y,z\in\N\}=\{2T_x+T_y+T_z:\ x,y,z\in\N\}=\N.$$
By [S15, Theorem 1.14], $T_x+T_y+2p_5(z)$ with $p_5(z)=z(3z-1)/2$ is also universal over $\Z$.
These, together with (1.1) and (1.2), indicate that (1.3) holds
for $(a,b,c)=(1,1,2),(1,2,2),(2,2,2),(2,2,3).$

In Section 2 we will prove Theorem 1.1 as well as the following related result.

\proclaim{Theorem 1.2} $(1.3)$ holds if $(a,b,c)$ is among
the following $7$ triples:
$$(1,2,3),\ (1,2,4),\ (1,2,5),\ (2,2,4),\ (2,2,5),\ (2,3,3),\ (2,3,4).$$
\endproclaim

In view of Theorems 1.1-1.2 and Remark 1.1, we have reduced the converse of Theorem 1.1 to our following conjecture.

\proclaim{Conjecture 1.1} $(1.3)$ holds if $(a,b,c)$ is among the following six triples:
$$(2,2,6),\ (2,3,5),\ (2,3,7),\ (2,3,8),\ (2,3,9),\ (2,3,10).$$
\endproclaim
\Remark\ 1.2. It is easy to show that $(1.3)$ holds for $(a,b,c)=(2,3,7)$ if and only if for any $n\in\N$
we can write $168n+41$ as $21x^2+14y^2+6z^2$ with $x,y,z\in\Z$.
\medskip

Inspired by (1.2), we want to know for what $a,b,c,d\in\N$ with $b\ls c\ls d\ls a$ we have
$$\{x(ax+b)+y(ay+c)+z(az+d):\ x,y,z\in\Z\}=\N.\tag1.5$$
We achieve this in the following theorem which will be proved in Section 3.

\proclaim{Theorem 1.3} Let $a>2$ be an integer and let $b,c,d\in\N$ with $b\ls c\ls d\ls a$.
Then $(1.5)$ holds if and only if $(a,b,c,d)$ is among the following five quadruples:
$$(3,0,1,2),\ (3,1,1,2),\ (3,1,2,2),\ (3,1,2,3),\ (4,1,2,3).\tag1.6$$
\endproclaim
\Remark\ 1.3. For $a\in\{1,2\}$ and $b,c,d\in\N$ with $b\ls c\ls d\ls a$, we can easily show that if (1.5) holds then
$(a,b,c,d)$ is among the following five quadruples:
$$(1,0,0,1),\ (1,0,1,1),\ (2,0,0,1),\ (2,0,1,1),\ (2,1,1,1).$$
The converse also holds since
$$x^2+y^2+2T_z,\ x^2+2T_y+2T_z,\ 2x^2+2y^2+T_z,\ 2x^2+T_y+T_z,\ T_x+T_y+T_z$$
are all universal over $\Z$ (cf. [S07]).
\medskip

We also note some other universal sums. For example, we have
$$\{x^2+y(3y+1)+z(3z+2): x,y,z\in\Z\}=\{x^2+y(4y+1)+z(4z+3): x,y,z\in\Z\}=\N$$
which can be easily proved.

Based on our computation, we formulate the following conjecture for further research.

\proclaim{Conjecture 1.2} {\rm (i)} Any positive integer $n\not=225$ can be written as $p(p-1)/2+q(q-1)/2+r(r-1)/2$ with $p$ prime and $q,r\in\Z^+$.

{\rm (ii)} Each $n\in\N$ can be written as $x^2+y(3y+1)/2+z(2z-1)$ with $x,y,z\in\N$. Also, any $n\in\N$ can be written as
 $x^2+y(3y+1)/2+z(5z+3)/2$ with $x,y,z\in\N$.

{\rm (iii)} Every $n\in\Z^+$ can be written as $x^3+y^2+T_z$ with $x,y\in\N$ and $z\in\Z^+$. We also have $\{x^2+y(y+1)+z(z^2+1):\ x,y,z\in\N\}=\N$.

{\rm (iv)} Any $n\in\N$ can be written as $x^4+y(3y+1)/2+z(7z+1)/2$ with $x,y,z\in\Z$.
\endproclaim
\medskip

\heading{2. Proofs of Theorems 1.1-1.2}\endheading

\medskip
\noindent {\it Proof of Theorem 1.1}. For $x\in\Z\sm\{0\}$, clearly $ax^2+x\gs |x|(a|x|-1)\gs a-1$.
As $1=x(ax+1)+y(by+1)+z(cz+1)$ for some $x,y,z\in\Z$, we must have $a\ls 2$.

{\it Case} 1. $a=b=1$.

As $1\not\in\{x(x+1)+y(y+1):\, x,y\in\Z\}$, we must have $1\in\{z(cz+1):\, z\in\Z\}$ and hence $c=2$.
(Note that if $c>2$ then $cz^2+z\gs c-1>1$ for all $z\in\Z\sm\{0\}$.)

{\it Case} 2. $a=1<b$.

If $b>2$, then $y(by+1)\gs b-1>1$ and $z(cz+1)\gs c-1>1$ for all $y,z\in\Z\sm\{0\}$.
As $1=x(x+1)+y(by+1)+z(cz+1)$ for some $x,y,z\in\Z$, we must have $b=2$.
It is easy to see that $4\not\in\{x(x+1)+y(2y+1):\,x,y\in\Z\}$. If $c>5$, then
$z(cz+1)\gs c-1>4$ for all $z\in\Z\sm\{0\}$. As $4=x(x+1)+y(2y+1)+z(cz+1)$ for some $x,y,z
\in\Z$, we must have $c\in\{2,3,4,5\}$.

{\it Case} 3. $a=b=2$.

In view of (1.1),
$$5\not\in\{T_x+T_y:\ x,y\in\N\}=\{x(2x+1)+y(2y+1):\ x,y\in\Z\}.$$
If $c>6$, then $z(cz+1)\gs c-1>5$ for all $z\in\Z\sm\{0\}$.
As $5=x(2x+1)+y(2y+1)+z(cz+1)$ for some $x,y,z
\in\Z$, we must have $c\in\{2,3,4,5,6\}$.

{\it Case} 4. $a=2<b$.

Clearly, $2\not\in\{x(2x+1):\ x\in\Z\}$.
If $b>3$, then $y(by+1)\gs b-1>2$ and $z(cz+1)\gs c-1>2$ for all $y,z\in\Z\sm\{0\}$.
As $2=x(2x+1)+y(by+1)+z(cz+1)$ for some $x,y,z
\in\Z$, we must have $b=3$. Note that $x(2x+1)+y(3y+1)\not=9$
for all $x,y\in\Z$. If $c>10$, then $z(cz+1)\gs c-1>9$ for all $z\in\Z\sm\{0\}$.
Since $9=x(2x+1)+y(3y+1)+z(cz+1)$ for some $x,y,z\in\Z$, we must have $c\ls 10$.
Note that $48\not=x(2x+1)+y(3y+1)+z(6z+1)$ for all $x,y,z\in\Z$. So $c\in\{3,4,5,7,8,9,10\}$.

In view of the above, we have completed the proof of Theorem 1.1. \qed

\proclaim{Lemma 2.1} Let $u$ and $v$ be integers with $u^2+v^2$ a positive multiple of $5$.
Then $u^2+v^2=x^2+y^2$ for some $x,y\in\Z$ with $5\nmid xy$.
\endproclaim
\Proof. Let $a$ be the $5$-adic order of $\gcd (u,v)$, and write $u=5^au_0$ and $v=5^av_0$ with $u_0,v_0\in\Z$ not all divisible by $5$.
Choose $\da,\ve\in\{\pm1\}$ such that $u_0'\not\eq2v_0'\pmod 5$, where $u_0'=\da u_0$ and $v_0'=\ve v_0$.
Clearly, $5^2(u_0^2+v_0^2)=u_1^2+v_1^2$, where $u_1=3u_0'+4v_0'$ and $v_1=4u_0'-3v_0'$. Note that $u_1$ and $v_1$ are not all divisible by $5$
since $u_1\not\eq v_1\pmod 5$. Continue this process, we finally write $u^2+v^2=5^{2a}(u_0^2+v_0^2)$ in the form $x^2+y^2$ with $x,y\in\Z$
not all divisible by $5$. As $x^2+y^2=u^2+v^2\eq0\pmod 5$, we must have $5\nmid xy$.
This concludes the proof. \qed

\medskip

With the help of Lemma 2.1, we are able to deduce the following result.

\proclaim{Lemma 2.2} For any $n\in\N$ and $r\in\{6,14\}$, we can write $20n+r$ as $5x^2+5y^2+z^2$ with $x,y,z\in\Z$ and $2\nmid z$.
\endproclaim
\Proof. As $20n+r\eq r\eq2\pmod 4$, by the Gauss-Legendre theorem we can write $20n+r$ as $(2w)^2+u^2+v^2$ with $u,v,w\in\Z$ and $2\nmid uv$.
If $(2w)^2\eq-r\pmod 5$, then $u^2+v^2\eq2r\pmod 5$ and hence $u^2\eq v^2\eq r\pmod 5$.
If $(2w)^2\eq r\pmod 5$, then $u^2+v^2\eq2\pmod 4$ is a positive multiple of $5$ and hence by Lemma 2.1 we can write it as $s^2+t^2$, where $s$ and $t$
are odd integers with $s^2\eq -r\pmod5$ and $t^2\eq r\pmod 5$. If $5\mid w$, then one of $u^2$ and $v^2$ is divisible by $5$ and the other is congruent to $r$ modulo $5$.

By the above, we can always write $20n+r=x^2+y^2+z^2$ with $x,y,z\in\Z$, $2\nmid z$ and $z^2\eq r\pmod 5$. Note that $x^2\eq-y^2=(\pm2y)^2\pmod 5$. Without loss of generality, we assume that $x\eq2y\pmod 5$ and hence $2x\eq-y\pmod5$.
Set $\bar x=(x-2y)/5$ and $\bar y=(2x+y)/5$. Then
$$20n+r=x^2+y^2+z^2=5\bar x^2+5\bar y^2+z^2.$$
This concludes the proof. \qed
\medskip

\Remark\ 2.1. Let $n\in\N$ and $r\in\{6,14\}$. In contrast with Lemma 2.2, we conjecture that $20n+r$ can be written as $5x^2+5y^2+(2z)^2$ with $x,y,z\in\Z$ unless
$r=6$ and $n\in\{0,11\}$, or $r=14$ and $n\in\{1,10\}$.

\proclaim{Lemma 2.3} {\rm (i)} For any positive integer $w=x^2+2y^2$ with $x,y\in\Z$, we can write $w$ in the form $u^2+2v^2$ with $u,v\in\Z$ such that $u$ or $v$ is not divisible by $3$.

{\rm (ii)} $w\in\N$ can be written as $3x^2+6y^2$ with $x,y\in\Z$, if and only if $3\mid w$ and $w=u^2+2v^2$ for some $u,v\in\Z$.

{\rm (iii)} Let $n\in\N$ with $6n+1$ not a square. Then, for any $\da\in\{0,1\}$ we can write $6n+1$ as $x^2+3y^2+6z^2$ with $x,y,z\in\Z$ and $x\eq\da\pmod2$.
\endproclaim
\Remark\ 2.2. Part (i) first appeared in the middle of a proof given on page 173 of [JP] (see also [S15, Lemma 2.1] for other similar results).
Parts (ii) and (iii) are Lemmas 3.1 and 3.3 of the author [S15].

\medskip\noindent{\it Proof of Theorem 1.2}. Let us fix a nonnegative integer $n$.

(i) As $24n+11\eq3\pmod 8$, by the Gauss-Legendre theorem there are odd integers $u,v,w$ such that $24n+11=u^2+v^2+w^2=w^2+2\bar u^2+2\bar v^2$,
where $\bar u=(u+v)/2$ and $\bar v=(u-v)/2$.
As $2(\bar u^2+\bar v^2)\eq 11-w^2\eq10\eq2\pmod 8$, we have $\bar u\not\eq\bar v\pmod2$.
Without loss of generality, we assume that $2\mid \bar u$ and $2\nmid \bar v$.
If $3\nmid\bar v$, then $\gcd(6,\bar v)=1$.
When $3\mid\bar v$, we have $3\nmid \bar u$ (since $w^2\not\eq11\pmod3$),
and $w^2+2\bar v^2$ is a positive multiple of $3$,
thus by Lemma 2.3(i) there are $s,t\in\Z$ with $3\nmid st$ such that $s^2+2t^2=w^2+2\bar v^2\eq3\pmod 8$ and hence $2\nmid st$.
Anyway, $24n+11$ can be written as $r^2+2s^2+2t^2$ with $r,s,t\in\Z$ and $\gcd(6,t)=1$. Since $r^2+2s^2\eq11-2t^2\eq0\pmod 3$,
by Lemma 2.3(ii) we may write $r^2+2s^2=3r_0^2+6s_0^2$ with $r_0,s_0\in\Z$. Since $3r_0^2+6s_0^2=r^2+2s^2\eq 11-2t^2\eq 9\pmod 8$,
we have $r_0^2+2s_0^2\eq3\pmod 8$ and hence $2\nmid r_0s_0$. Write $s_0=2x+1$, $r_0$ or $-r_0$ as $4y+1$, and $t$ or $-t$ as $6z+1$, where
$x,y,z\in\Z$. Then
$$24n+11=6(2x+1)^2+3(4y+1)^2+2(6z+1)^2$$
and hence $n=x(x+1)+y(2y+1)+z(3z+1)$. This proves (1.3) for $(a,b,c)=(1,2,3)$.

(ii) By the Gauss-Legendre theorem, there are $s,t,v\in\Z$ such that $32n+14=(2s+1)^2+(2t+1)^2+(2v)^2$
and hence $16n+7=(s+t+1)^2+(s-t)^2+2v^2$. As one of $s+t+1$ and $s-t$ is even, we have $16n+7=(2u)^2+w^2+2v^2$
for some $u,w\in\Z$. Clearly $2\nmid w$, $2v^2\eq 7-w^2\eq2\pmod 4$,
 and $4u^2\eq 7-2v^2-w^2\eq4\pmod 8$. So, $u,v,w$ are all odd. Note that $w^2\eq 7-4u^2-2v^2\eq 7-4-2=1\pmod{16}$
 and hence $w\eq\pm1\pmod 8$. Now we can write $u$ as $2x+1$, $v$ or $-v$ as $4y+1$, $w$ or $-w$ as $8z+1$, where
 $x,y,z$ are integers. Thus
 $$16n+7=4(2x+1)^2+2(4y+1)^2+(8z+1)^2$$
 and hence $n=x(x+1)+y(2y+1)+z(4z+1)$. This proves (1.3) for $(a,b,c)=(1,2,4)$.

(iii) By Dickson [D39, pp.\,112-113] (or [JKS]),
$$E(10x^2+5y^2+2z^2)=\{8q+3:\ q\in\N\}\cup\bigcup_{k,l\in\N}\{25^k(5l+1),25^k(5l+4)\}.$$
So, there are $u,v,w\in\Z$ such that $40n+17=10u^2+5v^2+2w^2$. Clearly, $2\nmid v$, $2u^2+2w^2\eq 17-5v^2\eq4\pmod 8$
and hence $2\nmid uw$. Note that $2w^2\eq 17\eq2\pmod 5$ and hence $w\eq\pm1\pmod 5$. Thus, we can write
$u=2x+1$, $v$ or $-v$ as $4y+1$, and $w$ or $-w$ as $10z+1$, where $x,y,z$ are integers. Now we have
$$40n+17=10(2x+1)^2+5(4y+1)^2+2(10z+1)^2$$
and hence $n=x(x+1)+y(2y+1)+z(5z+1)$. This proves (1.3) for $(a,b,c)=(1,2,5)$.

(iv) By the Gauss-Legendre theorem, there are $u,v,w\in\Z$ with $2\nmid w$ such that
$$16n+5=(2u)^2+(2v)^2+w^2=2(u+v)^2+2(u-v)^2+w^2.$$
As $w^2\eq1\not\eq5\pmod 8$, both $u+v$ and $u-v$ are odd.
Since $w^2\eq5-2-2=1\pmod{16}$, we have $w\eq\pm1\pmod 8$.
Now we can write $u+v$ or $-u-v$ as $4x+1$, $u-v$ or $v-u$ as $4y+1$, and $w$ or $-w$ as $8z+1$, where $x,y,z\in\Z$.
Thus
$$16n+5=2(4x+1)^2+2(4y+1)^2+(8z+1)^2$$
and hence $n=x(2x+1)+y(2y+1)+z(4z+1)$.
This proves (1.3) for $(a,b,c)=(2,2,4)$.

(v) By Lemma 2.2, there are $u,v,w\in\Z$ with $2\nmid w$ such that $20n+6=5u^2+5v^2+w^2$.
Clearly, $u\not\eq v\pmod 2$, $w^2\eq 1\pmod 5$ and hence $w\eq\pm1\pmod 5$. Thus $w$ or $-w$
has the form $10z+1$ with $z\in\Z$. Observe that
$$40n+12=10u^2+10v^2+2w^2=5(u+v)^2+5(u-v)^2+2(10z+1)^2.$$
As $u+v$ and $u-v$ are both odd, we may write $u+v$ or $-u-v$ as $4x+1$, and $u-v$ or $v-u$ as $4y+1$,
where $x$ and $y$ are integers. Then
$$40n+12=5(4x+1)^2+5(4y+1)^2+2(10z+1)^2$$
and hence $n=x(2x+1)+y(2y+1)+z(5z+1)$. This proves (1.3) for $(a,b,c)=(2,2,5)$.

(vi) By Dickson [D39, pp.\,112-113],
$$E(x^2+y^2+3z^2)=\{9^k(9l+6):\ k,l\in\N\}.$$
So there are $u,v,w\in\Z$ such that $24n+7=u^2+v^2+3w^2$. As $u^2+v^2\not\eq 7\pmod 4$, we have $2\nmid w$
and hence $s=(u+v)/2\in\Z$ and $t=(u-v)/2\in\Z$. Now $24n+7=2s^2+2t^2+3w^2$. As $2(s^2+t^2)\eq 7-3w^2\eq 4\pmod 8$,
we have $s^2+t^2\eq2\pmod 4$ and hence $2\nmid st$. Note that $s^2+t^2\eq (7-3)/2=2\pmod 3$ and hence $3\nmid st$.
Now we can write $w$ or $-w$ as $4x+1$, $s$ or $-s$ as $6y+1$, $t$ or $-t$ as $6z+1$, where $x,y,z$ are integers.
Then
$$24n+7=3(4x+1)^2+2(6y+1)^2+2(6z+1)^2$$
and hence
$n=x(2x+1)+y(3y+1)+z(3z+1)$. This proves (1.3) for $(a,b,c)=(2,3,3)$.

(vii) Note that $48n+13\eq1\pmod 6$ but $48n+13\not\eq1\pmod 8$. By Lemma 2.3(iii), there are $u,v,w\in\Z$ such that
$48n+13=6u^2+(2v)^2+3w^2$. Clearly, $2\nmid w$ and $3\nmid v$. As $6u^2\eq 13-3\eq6\pmod 4$, we must have $2\nmid u$.
Since $4v^2\eq 13-6u^2-3w^2\eq 4\pmod 8$, we have $2\nmid v$. Observe that
$$3w^2\eq 13-6u^2-4v^2\eq 13-6-4=3\pmod{16}$$
and hence $w\eq\pm1\pmod8$. Now we can write $u$ or $-u$ as $4x+1$,
$v$ or $-v$ as $6y+1$, and $w$ or $-w$ as $8z+1$, where $x,y,z\in\Z$.
Thus
$$48n+13=6(4x+1)^2+4(6y+1)^2+3(8z+1)^2$$
and hence $n=x(2x+1)+y(3y+1)+z(4z+1)$. This proves (1.3) for $(a,b,c)=(2,3,4)$.

So far we have completed the proof of Theorem 1.2. \qed

\heading{3. Proof of Theorems 1.3}\endheading

\proclaim{Lemma 3.1} For any positive integer $n$, we can write $6n+1$ as $x^2+y^2+2z^2$ with $x,y,z\in\Z$ and $3\nmid xyz$.
\endproclaim
\Remark\ 3.1. This is [S16, Lemma 4.3(ii)] proved by the author with the help of a result in [CL]. Combining it with Lemma 2.3(ii), for any
$n\in\Z^+$ and $\da\in\{0,1\}$ we can write $6n+1$ as $x^2+3y^2+6z^2$ with $x,y,z\in\Z$ and $x\eq\delta\pmod2$, which extends Lemma 2.3(iii) and confirms a conjecture 
in [S15, Remark 3.4].

\medskip\noindent{\it Proof of Theorem 1.3}. (i) If $|x|\gs2$, then
$$x(ax+b)\gs |x|(a|x|-b)\gs 2(2a-b)\gs 2a,$$
and similarly $x(ax+c)\gs 2a$ and $x(ax+d)\gs 2a$.
So, if (1.5) holds then we must have
$$\{0,1,\ldots,2a-1\}\se\{x(ax+b)+y(ay+c)+z(az+d):\ x,y,z\in\{0,\pm1\}\}$$
and hence
$$2a\ls|\{x(ax+b)+y(ay+c)+z(az+d):\ x,y,z\in\{0,\pm1\}\}|\ls 3^3=27.$$
Note that $a\in\{3,4,\ldots,13\}$ and $0\ls b\ls c\ls d\ls a$. Via a computer we
find that if $(a,b,c,d)$ is not among the five quadruples in (1.6) then one of $1,2,\ldots,17$ cannot be written as
$x(ax+b)+y(ay+c)+z(az+d)$ with $x,y,z\in\Z$. For example, $x(4x+2)+y(4y+2)+z(4z+3)\not=17$ for any $x,y,z\in\Z$.
This proves the ``only if" part of Theorem 1.3.

(ii) Now we turn to prove the ``if" part of Theorem 1.3. Let us fix a nonnegative integer $n$.

(a) By [S15, Theorem 1.7(iv)], there are $u,v,x\in\Z$ such that $12n+5=u^2+v^2+36x^2$.
Clearly $u\not\eq v\pmod 2$ and $3\nmid uv$. Without loss of generality, we assume that
$u\eq\pm1\pmod6$ and $v\eq\pm2\pmod 6$. We may write $u$ or $-u$ as $6y+1$, and $v$ or $-v$ as $6z+2$, where $y$ and $z$
are integers. Thus
$$12n+5=36x^2+(6y+1)^2+(6z+2)^2$$
and hence $n=3x^2+y(3y+1)+z(3z+2)$. This proves (1.5) for $(a,b,c,d)=(3,0,1,2)$.

(b) Let $\da\in\{0,1\}$. By the Gauss-Legendre theorem, $12n+6+3\da$ can be written as the sum of three squares.
In view of [S16, Lemma 2.2], there are $u,v,w\in\Z$ with $3\nmid uvw$ such that $12n+6+3\da=u^2+v^2+w^2$.
Clearly, $u,v,w$ are neither all odd nor all even. Without loss of generality, we assume that $2\nmid u$ and $2\mid w$.
Then $v\not\eq \da\pmod2$. Obviously, $u\eq\pm1\pmod 6$, $v\eq\pm(1+\da)\pmod 6$ and $w\eq\pm2\pmod 6$. Thus we may write $u$ or $-u$ as $6x+1$,
$v$ or $-v$ as $6y+1+\da$, and $w$ or $-w$ as $6z+2$, where $x,y,z\in\Z$. Therefore,
$$12n+6+3\da=(6x+1)^2+(6y+1+\da)^2+(6z+2)^2$$
and hence $n=x(3x+1)+y(3y+1+\da)+z(3z+2)$. This proves (1.5) for $(a,b,c,d)=(3,1,1,2),(3,1,2,2)$.

(c) By Lemma 3.1, there are $u,v,w\in\Z$ with $3\nmid uvw$ such that $6n+7=u^2+v^2+2w^2$ and hence
$12n+14=(u+v)^2+(u-v)^2+(2w)^2$. As $(u+v)^2+(u-v)^2\eq2\pmod 4$, both $u+v$ and $u-v$ are odd.
Since $(u+v)^2+(u-v)^2\eq 14-1\eq1\pmod3$, without loss of generality we may assume that $u+v\eq\pm1\pmod6$ and $u-v\eq3\pmod 6$.
Now we may write $u+v$ or $-u-v$ as $6x+1$, $w$ or $-w$ as $3y+1$, and $u-v$ as $6z+3$, where $x,y,z$ are integers.
Then
$$12n+14=(6x+1)^2+(6y+2)^2+(6z+3)^2$$
and hence $n=x(3x+1)+y(3y+2)+z(3z+3)$. This proves (1.5) for $(a,b,c,d)=(3,1,2,3)$.

(d) As $16n+14\eq2\pmod 4$, by the Gauss-Legendre theorem $16n+14=u^2+v^2+w^2$ for some $u,v,w\in\Z$ with $2\nmid uv$ and $2\mid w$.
Since $w^2\eq 14-u^2-v^2\eq12\eq4\pmod 8$, $w/2$ or $-w/2$ has the form $4y+1$ with $y\in\Z$. Thus
$u^2+v^2\eq 14-4(w/2)^2\eq10\pmod{16}$. Without loss of generality, we assume that $u\eq\pm1\pmod 8$ and $v\eq\pm3\pmod 8$.
Write $u$ or $-u$ as $8x+1$, and $v$ or $-v$ as $8z+3$, where $x,z\in\Z$. Then
$$16n+14=(8x+1)^2+(8y+2)^2+(8z+3)^2$$
and hence $n=x(4x+1)+y(4y+2)+z(4z+3)$.  This proves (1.5) for $(a,b,c,d)=(4,1,2,3)$.

In view of the above, we have completed the proof of Theorem 1.3. \qed

\medskip
\Ack. The author would like to thank Dr. Hao Pan for his comments on the proof of Lemma 2.2, and the referee for his/her helpful suggestions.

\enddocument